\documentclass{amsart}

\usepackage{amsfonts}
\usepackage{cases}
\usepackage{mathrsfs}
\usepackage{bbm}
\usepackage{amssymb}
\usepackage{txfonts}
\usepackage{amscd}
\usepackage{amsfonts,latexsym,amsmath,amsthm,amsxtra,mathdots}
\usepackage[bookmarks,colorlinks]{hyperref}
\usepackage[all,cmtip]{xy}
\RequirePackage{amsmath} \RequirePackage{amssymb}
\usepackage{color}
\usepackage{colordvi}
\usepackage{multicol}
\usepackage{tikz}
\usepackage{hyperref}
\usepackage{graphicx}
\usepackage{amsmath}
\usepackage{amsmath,amscd}
\usetikzlibrary{matrix,arrows}
\usepackage{textcomp}

\makeatletter
\DeclareFontFamily{OMX}{MnSymbolE}{}
\DeclareSymbolFont{MnLargeSymbols}{OMX}{MnSymbolE}{m}{n}
\SetSymbolFont{MnLargeSymbols}{bold}{OMX}{MnSymbolE}{b}{n}
\DeclareFontShape{OMX}{MnSymbolE}{m}{n}{
    <-6>  MnSymbolE5
   <6-7>  MnSymbolE6
   <7-8>  MnSymbolE7
   <8-9>  MnSymbolE8
   <9-10> MnSymbolE9
  <10-12> MnSymbolE10
  <12->   MnSymbolE12
}{}
\DeclareFontShape{OMX}{MnSymbolE}{b}{n}{
    <-6>  MnSymbolE-Bold5
   <6-7>  MnSymbolE-Bold6
   <7-8>  MnSymbolE-Bold7
   <8-9>  MnSymbolE-Bold8
   <9-10> MnSymbolE-Bold9
  <10-12> MnSymbolE-Bold10
  <12->   MnSymbolE-Bold12
}{}
\let\llangle\@undefined
\let\rrangle\@undefined
\DeclareMathDelimiter{\llangle}{\mathopen}%
                     {MnLargeSymbols}{'164}{MnLargeSymbols}{'164}
\DeclareMathDelimiter{\rrangle}{\mathclose}%
                     {MnLargeSymbols}{'171}{MnLargeSymbols}{'171}
\makeatother

\newcommand{\lk}{\operatorname{lk}}

\hypersetup{
    colorlinks,
   linkcolor={red!10!black},
   citecolor={blue!10!black},
    urlcolor={blue!10!black}
}

     \newcommand{\BZ}{{\mathbb {Z}}}

     \newcommand{\Aut}{{\mathrm{Aut}}}

\def\-{^{-1}}

\newcommand{\delete}[1]{}

    \newcommand{\st}{{\mathrm{st}}}

    \theoremstyle{plain}

\newtheorem{thm}{Theorem}[section]
\newtheorem{lem}[thm]{Lemma}
\newtheorem{prop}[thm]{Proposition}
\newtheorem{cor}[thm]{Corollary}

\newtheorem*{thmA}{Theorem A}
\newtheorem*{corB}{Corollary B}

\newtheorem*{thmC}{Theorem C}
\newtheorem*{rem1}{Remark 1}
\newtheorem*{rem2}{Remark 2}
\newtheorem*{rem3}{Remark 3}

    \numberwithin{equation}{section}

\def\Proof{\noindent{\bf Proof}\quad}
\def\qed{\hfill$\square$\smallskip}

\begin{document}

\title{The Farrell--Jones Conjecture for normally poly-free groups}

\author{Benjamin Br\"uck}
\address{Fakult\"at f\"ur Mathematik, Universit\"at Bielefeld, Postfach 100131, D-33501 Bielefeld, Germany}
\email{benjamin.brueck.maths@gmail.com}

\author{Dawid Kielak}
\address{Fakult\"at f\"ur Mathematik, Universit\"at Bielefeld, Postfach 100131, D-33501 Bielefeld, Germany}
\email{dkielak@math.uni-bielefeld.de}

\author{Xiaolei Wu}
\address{ Mathematische Institut, Universit\"at Bonn, Endenicher Allee 60, D-53115, Bonn, Germany}
\email{xwu@math.uni-bonn.de}


\date{September, 2020}


\begin{abstract}
We 
prove the $K$- and $L$-theoretic Farrell--Jones Conjecture with coefficients in an additive category for every  normally poly-free group, in particular for 
even Artin groups of FC-type, and for all groups of the form $A\rtimes \BZ$ where $A$ is a right-angled Artin group. Our proof relies on the work of Bestvina--Fujiwara--Wigglesworth on the Farrell--Jones Conjecture for free-by-cyclic groups.
\end{abstract}

\maketitle

\section*{Introduction}

Recently, Bestvina, Fujiwara, and Wigglesworth proved
the Farrell--Jones Conjecture for free-by-cyclic groups where the free group has finite rank \cite{BFW19}. The purpose of this note is to extend their result and to allow the free group to have infinite rank. More generally, we will deal with all normally poly-free groups. Recall that a group $G$ is called \emph{normally poly-free} if it admits a chain of subgroups $ 1 =G_0\leq G_1\leq \cdots \leq G_{n-1}\leq G_n=G$, such that each $G_i$ is  normal in  $G$ and $G_i/G_{i-1}$ is a free group of possibly infinite rank.
Our main theorem can be stated as follows.

\begin{thmA} \label{thm}
The K- and L-theoretic Farrell--Jones Conjecture   with finite wreath products and coefficients in an additive category is true for normally poly-free groups.
\end{thmA}

For more information about the Farrell--Jones Conjecture with finite wreath products and coefficients in an additive category (which we will abbreviate to FJCw) we refer to \cite[Section 2.3]{We15}.

Normally poly-free groups have two related families. The more classical one is the family of \emph{poly-free} groups, where each $G_i$ is only required to be normal in $G_{i+1}$. 
The second one, defined by  Aravinda, Farrell, and Roushon \cite{AFR00}, consists of \emph{strongly poly-free} groups, which are normally poly-free and where each quotient $G_i/G_{i-1}$ is a free group of finite rank, such that the conjugation action of $G/G_{i-1}$ on $G_i/G_{i-1}$ can be realized by some surface homeomorphism. Aravinda, Farrell, and Roushon showed that FJCw holds for  all strongly poly-free groups. As a corollary, they proved the conjecture for pure braid groups. Thus one might expect that our theorem can be used to verify the Farrell--Jones Conjecture for more Artin groups. Indeed, in \cite{BMP19}, Blasco-Garcia, Mart\'inez-P\'erez, and  Paris showed that even Artin groups of FC-type (see subsection \ref{artin-group} for the definition) are poly-free. A careful reading of their paper, in particular \cite[Proposition 3.2]{BMP19},  shows  that even Artin groups of FC-type are actually normally poly-free. Hence we have the following.

\begin{corB}
The K- and L-theoretic Farrell--Jones Conjecture   with finite wreath products and coefficients in an additive category is true for every even Artin group of FC-type.
\end{corB}

\begin{rem1}
Osajda and  Huang have independently obtained a proof of the Farrell--Jones Conjecture for all Artin groups of FC-type \cite{HO19}. Sayed K. Roushon has also proved the conjecture for some classes of Artin groups \cite{Ro18}. In \cite{Wu19}, the third author proved the conjecture for more even Artin groups using Theorem A. 
\end{rem1}

\begin{rem2}
Using \cite[Theorem 8.7 (2)]{Lu18}, one can prove by induction that every poly-free group satisfies the Baum-Connes Conjecture with coefficients. In particular, every even Artin group of FC-type satisfies the Baum-Connes Conjecture with coefficients.
\end{rem2}

Potentially, our theorem could be applied to more Artin groups. In \cite[Question 2]{Be99} and the discussions below it, Bestvina asks whether all Artin groups are virtually poly-free. So one might further ask whether all Artin groups are virtually normally poly-free. To the best of the authors' knowledge, this question remains open.

Our main theorem offers another application. Given a simplicial graph $\Gamma$, we can associate to it a right-angled Artin group $A_\Gamma$ defined as follows: the generators of $A_\Gamma$ are the vertices of $\Gamma$; for any two vertices $v$ and $w$ which are connected by an edge, we add a relation $[v,w]=1$. Note that right-angled Artin groups are $\mathrm{CAT}(0)$ groups, so, in particular, they satisfies FJCw \cite{BL12, We12}. With the help of Theorem A, we are able to extend this to the following.

\begin{thmC}\label{thm:FJCw-raag-by-cyclic}
Let $\Gamma$ be a finite simplicial graph and let $A_\Gamma$ be the corresponding  right-angled Artin group. The K- and L-theoretic Farrell--Jones Conjecture   with finite wreath products and coefficients in an additive category is true for every group of the form $A_\Gamma \rtimes_f\BZ$, where $f$ is an automorphism of the group $A_\Gamma$.
\end{thmC}

\begin{rem3}
Note that our results also hold for the Farrell--Jones Conjecture in A-theory since all the results and inheritance properties we used in our proof also holds for A-theory \cite{ELM+18,FW17, KUW18}.
\end{rem3}


\textbf{Acknowledgements.} 
Br\"uck and Kielak were supported by grants BU 1224/2-1 and KI 1853/3-1, respectively,  within the Priority Programme 2026 `Geometry at infinity' of the German Science Foundation (DFG).

Wu was partially supported by Wolfgang L\"uck's ERC Advanced Grant “KL2MG-interactions”
(no. 662400) and the DFG Grant under Germany's Excellence Strategy - GZ 2047/1, Projekt-ID 390685813. Wu wants to thank Grigori Avramidi, Daniel Kasprowski and Henrik R\"uping for some stimulating discussions. He also wants to thank Daniel Kasprowski, Wolfgang L\"uck and Christoph Winges for some helpful comments on previous versions of the paper. Part of this work was done when Wu was visiting Bielefeld; he wants to thank Kai-Uwe Bux and his group for the warm hospitality.

\section{Preliminaries}

\subsection{Farrell--Jones Conjecture.} We list some inheritance properties and results on FJCw which we need to prove our main results; these results can be found for example in \cite[Section 2.3]{We15} and \cite{GMR15}, see also \cite{RV18} for a survey about the current status of the conjecture. The last item is the recent result of Bestvina, Fujiwara and Wigglesworth. Note that they did not state their result for the Farrell-Jones Conjecture with finite wreath product version. But since their proof only relies  on inheritance properties, Bartels' result on relative hyperbolic groups  and Knopf's result on groups acting acylindrically  on trees, by \cite[Remark 4.7]{Ba17} and \cite[Corollary 4.4]{Kn19}, their result also holds for FJCw. 
\begin{thm}\label{qut}
(1) If a group $G$ satisfies FJCw, then every subgroup $H_1 \leq G$ and every finite index overgroup $H_2 \geq G$ satisfies FJCw.\\
(2) If $G_1$ and $G_2$ satisfy FJCw, then the direct product $G_1 \times G_2$ and the free product $G_1 \ast G_2$ satisfy FJCw.\\
(3) Let $\{G_i \mid  i \in I\}$ be a directed system of groups (with not necessarily injective
structure maps). If each $G_i$ satisfies FJCw, then the colimit $colim_{i \in I} G_i$
satisfies FJCw.\\
(4) Let $\phi \colon G \rightarrow Q$ be a
group homomorphism. If $Q$, $Ker(\phi)$ and $\phi^{-1}(C)$ satisfy FJCw for every infinite cyclic subgroup $C \leq Q$, then G satisfies FJCw.\\
(5) CAT(0) groups satisfy FJCw.\\
(6) Virtually solvable groups satisfy FJCw.\\
(7) Fundamental groups of graphs of abelian groups satisfy FJCw.\\
(8) $F\rtimes \BZ$ satisfies FJCw, where $F$ is a free group of finite rank.
\end{thm}

Note that in (7) above we are not assuming the underlying graph to be finite.

\subsection{Artin groups.}\label{artin-group} We give a quick review of Artin groups, see \cite{Mc} for more information. An \textbf{Artin group} (or generalized braid group) $A$ is a group with a presentation of the form
\[\langle x_1,x_2,\cdots, x_n \mid \langle x_1,x_2\rangle^{m_{1,2}} =\langle x_2,x_1\rangle^{m_{2,1}},\cdots,\langle x_{n-1},x_{n}\rangle^{m_{n-1,n}} =\langle x_{n},x_{n-1}\rangle^{m_{n,n-1}} \rangle  \]
where $m_{i,j}=m_{j,i}\in\{2,3,\dots,\infty\},i<j$, and where for $m_{i,j}\in\{2,3,\dots\}$ the symbol $\langle x_i,x_j\rangle^{m_{i,j}}$ denotes an alternating product of $x_i$ and $x_j$ of length $m$, beginning with $x_i$. For example $\langle x_1,x_3\rangle^3 =x_1x_3x_1$.  When $m_{i,j}=\infty$, there is (by convention) no relation for $x_i$ and $x_j$.   An Artin group is called \textbf{even} if $m_{ij}$ is an even integer or $\infty$ for every $i,j$. It is a \textbf{right-angled Artin group} (RAAG) if $m_{ij}$ is either $2$ or $\infty$ for every $i,j$.

This data can be encoded by a \textbf{Coxeter diagram}: a labeled graph with $n$ vertices
$x_1,\cdots, x_n$, where two vertices $x_i$ and $x_j$ are connected by an edge if $2\leq m_{i,j}<  \infty$ and edges are labeled by $m_{i,j}$ whenever $m_{i,j} \geq 3$. We will adopt the convention used for right-angled Artin groups and leave out the label when it is $2$. If $m_{i,j} = \infty$, we will not connect the vertices by an edge.

Given an Artin group with the above presentation, one further obtains a Coxeter group $W$ by adding the relation $x_i^2 =1$ for all $i$. We say that the Artin group $A$ is of \textbf{spherical type} if
the associated Coxeter group $W$ is finite.

Given a Coxeter graph $\Gamma$ with associated Artin group $A$ and Coxeter group $W$, and a subset $T$ of $\{x_1,x_2,\cdots, x_n\}$, we denote by $A_T$ (resp.~$W_T$) the subgroup of $A$ (resp.~$W$)  generated by $T$, and by $\Gamma_T$ the full subgraph of $\Gamma$ spanned by $T$. Here each edge of $\Gamma_T$ is labeled with the same number as its corresponding edge of $\Gamma$. The group $W_T$ is  the Coxeter group of $\Gamma_T$ \cite{Bo68}, and $A_T$ is the Artin group of $\Gamma_T$ \cite{Va83}.  Now a subset $T$ of $\{x_1,x_2,\cdots, x_n\}$ is called \textbf{free of infinity} if $m_{s,t} \neq \infty$ for all $s,t\in T$. We say that the Artin group $A$ is of \textbf{FC-type} if $A_T$ is of spherical type for every free of infinity subset $T$ of the generating set $\{x_1,x_2,\cdots, x_n\}$.

\section{FJCw for normally poly-free groups} \label{Sec:npf-group}

In this section we prove that FJCw holds for normally poly-free groups.

\begin{lem} \label{HNN-auto-FJCw}
Let $G$ be a group and $f$ be an automorphism of $G$. If  the semidirect product $G\rtimes_{f} \BZ$ satisfies FJCw, then for every subgroup $H$ of $G$, the HNN-extension \[G_H = \langle G,t \mid tht^{-1}= f(h),h\in H \rangle\] also satisfies FJCw.

\end{lem}
\Proof 
Let $q \colon G_H \to G\rtimes_f \BZ$ denote the natural quotient map.
Since $G_H$ is an HNN-extension, it acts on the corresponding Bass--Serre tree $T$ with vertex stabilizers all conjugate to $G$. Now the kernel of $q$ as a subgroup of $G_H$ also acts on $T$. Moreover, it acts on $T$ with trivial vertex stabilizers since $q$ maps the vertex stabilizer $G$ (and every conjugate of $G$) injectively into $G\rtimes \BZ$. Hence the kernel of $q$ is a free group, which satisfies FJCw. 

Now given any infinite cyclic subgroup $C$ of $G\rtimes \BZ$, if $C$ does not lie in $G\rtimes 0$, then $q^{-1}(C)$ has trivial intersection with every conjugate of $G$, so it acts freely on $T$. Hence $q^{-1}(C)$ is also a free group in this case. On the other hand, if $q^{-1}(C)\leq G\rtimes 0$, then $q^{-1}(C)$ acts on $T$ with cyclic stabilizers since $q$ restricted to $G$ is injective (hence $q$ restricted to any vertex stabilizer is injective). Thus $q^{-1}(C)$ is the fundamental group of some graph of cyclic groups. By Theorem \ref{qut} (7), the group $q^{-1}(C)$ satisfies FJCw. Now since $G\rtimes_f\BZ$  also satisfies FJCw, Theorem \ref{qut} (4) implies that $G_H$ satisfies FJCw. \qed

Note that the same proof also shows the following proposition, which might be of independent interest.

\begin{prop}\label{quotient-HNN-FJC}
Let $G $ be an HNN-extension of $H$ along two subgroups. Suppose that  we have a map $q\colon G\rightarrow N$ such that $q$ is injective when restricted to $H$. If $N$ satisfies FJCw, then $G$ also satisfies FJCw.
\end{prop}

\begin{cor}
Let  $G$ be a group which satisfies FJCw. Then for every subgroup $H$ of $G$, the HNN-extension $\langle G, t \mid  tht^{-1}= h,h \in H \rangle$ also satisfies FJCw.
\end{cor}

\Proof This follows from the fact that $G\times \BZ$ also satisfes FJCw.
\qed

In order to extend the Bestvina--Fujiwara--Wigglesworth result from finite rank free-by-cyclic groups to the infinite rank case, we need the following theorem.

\begin{thm}\cite[Theorem 4.4]{FJ50} \label{thm-free-factor}
If $F$ is a free group, $F = A \ast B $, $H$ is a subgroup of $F$ and if
$A \leq H$, then there is a subgroup $C$ of $H$ for which $H = A \ast C$.
\end{thm}

We are now ready to prove the following.

\begin{thm}\label{fjc-fbc}
Let $F$ be a (countable) free group of  infinite rank and $f$ an automorphism of $F$. Then the group $G=F\rtimes_f \BZ$ satisfies FJCw.
\end{thm}
\Proof 
Assume  the free basis of $F$ is $\{x_1,x_2,\dots, x_i,\dots \}$ and denote the subgroup generated by $\{x_1,x_2,\dots, x_n\}$ by $F_n$. Given $i>0$, let
$n_i\geq i$ be the smallest natural number such that $f(x_j)$ is contained in the subgroup $F_{n_i}$ for every $1\leq j\leq i$. Let 
$$f_i\colon F_i \rightarrow F_{n_i}$$
be the restriction of $f$ to $F_i$. Now let $G_i$ be the corresponding HNN-extension of $F_{n_i}$, i.e.
$$G_i = \langle F_{n_i},t_i \mid t_ixt_i^{-1}= f_i(x),\text{ for any } x\in F_i \rangle$$
Notice that $G$ is a colimit of $G_i$. So by Theorem \ref{qut} (3), we only need to show that each $G_i$ satisfies FJCw.

Now let $\bar{F}_{i}$ be the free subgroup of $F$ generated by $\{ x_{i +1},x_{i+2},\ldots \}$. Then $F = F_{i} \ast \bar{F}_{i}$. Since $f$ is an automorphism of $F$, we also have $F =  f(F_{i}) \ast f(\bar{F}_{i})$. Now $F_{n_i}$ is a subgroup of $F$ and $f(F_i) \leq F_{n_i}$, and so by Theorem \ref{thm-free-factor} there is a subgroup $C_i$ of $F_{n_i}$ such that $F_{n_i} = f(F_i)\ast C_i$. Notice that $C_i$ is a free group of rank $n_i -i$. Choosing a free basis of $C_i$, we can now complete the injection map $f_i:F_i \rightarrow F_{n_i}$ to an automorphism of  $F_{n_i}$ by mapping the remaining basis elements $x_{i+1},\cdots, x_{n_i}$ bijectively to the  basis of $C_i$. Denoting the map by $\bar{f}_i$, we have a free-by-cyclic group $ F_{n_i}\rtimes_{\bar{f}_i} \BZ$ which in particular  satisfies FJCw. By Lemma  \ref{HNN-auto-FJCw}, the groups $G_i$ also satisfy FJCw for every $i$.
\qed

\textbf{Proof of Theorem A.} Recall a group $G$ is called normally poly-free if it has a chain of subgroups $1 =G_0\leq G_1\leq G_2\cdots G_{n-1}\leq G_n=G$, such that each $G_i$ is  normal in $G$ and $G_i/G_{i-1}$ is a free group of possibly infinite rank. We prove FJCw for $G$ by induction on the length $n$. Suppose FJCw is true for all normally poly-free groups of length $n-1$. We have the following short exact sequence for $G$:
$$\begin{CD}
0 @> >> G_1 @>i>>  G  @> q >>G/G_1 @>  >>0, 
\end{CD}$$

Now $G/G_1$ has length $n-1$ and satisfies FJCw by induction. Thus by Theorem \ref{qut} (4), we only need to show that for every infinite cyclic subgroup $C$ in $G/G_1$, the group $q^{-1}(C)$ satisfies FJCw. But since $G_1$ is a free group, the group $q^{-1}(C)$ is free-by-cyclic. And now by Theorem \ref{fjc-fbc}, $q^{-1}(C)$ satisfies FJCw.

\section{FJCw for RAAG-by-cyclic groups}\label{sec:FJCw-raag-by-cyc}
In this section we prove FJCw for RAAG-by-cyclic groups. We first need some discussions about RAAGs and their automorphisms and refer to \cite{CV09} for more details. 

Given a simplicial graph $\Gamma$, let $\lk(v)$ and $\st(v)$ denote the link and star of a vertex $v$ of $\Gamma$, respectively.
An important tool for studying automorphisms of RAAGs is the so-called \textbf{standard ordering} on the vertex set $V(\Gamma)$ of $\Gamma$. It is the partial pre-order that is given by $v\leq w$ if and only if $\lk(v) \subseteq \st(w)$.
This induces an equivalence relation on $V(\Gamma)$ whose equivalence classes are denoted by $[\cdot]$. The standard ordering induces a partial order on the equivalence classes where we say $[v]\leq [w]$ if $v\leq w$ (this does not depend on the choice of representatives).
It is easy to check that for all $v\in V(\Gamma)$, the subgroup $A_{[v]}$ is either a free or a free abelian group of rank equal to the size of $[v]$. Furthermore, $A_{\st[v]}\cong A_{[v]}\times A_{\lk (v)}$, where $\st[v]=\lk[v] \cup [v]$.
For more details about this ordering and the equivalence relation, see \cite[Section 2]{CV09}.

In the sequel, given a subgroup $B \leq A_\Gamma$ and an element $g \in A_\Gamma$, we will use $B^g$ to denote the conjugate of $B$ by $g$.

\begin{lem}\label{structure-normal-closure}
Let $T=[v]$ be a standard equivalence class of $V(\Gamma)$. The normal closure $\llangle A_T \rrangle$ of $A_T$ in $A_\Gamma$ is given by the free product $\llangle A_T \rrangle = \left(\ast_{g\in X} A_T^g\right) \ast F$, where $F$ is a free group and $X$ is a system of coset representatives of $\llangle A_T \rrangle \backslash A_\Gamma / A_{\st[v]}$
\end{lem}
\Proof
The group $A_\Gamma$ decomposes as an amalgamated product $A_\Gamma=A_{\st[v]} \ast_{A_{\lk[v]}} A_{V(\Gamma)\backslash T}$, where $\lk[v]=\bigcap_{v'\in [v]} \lk (v')$ and $\st[v]=\lk[v] \cup [v]$. 
Let $q \colon A_\Gamma \to A_{V(\Gamma)\backslash T}$ denote the natural quotient map obtained by killing the generators in $T$. Note that the obvious embedding $A_{V(\Gamma)\backslash T} \to A_\Gamma$ shows that $q$ is in fact a retract.
The kernel of $q$ is precisely the normal closure of $A_T$. Hence we have $\llangle A_T \rrangle \cap A_{V(\Gamma)\backslash T} = \{ 1\}$. In particular, $\llangle A_T \rrangle \cap A_{\lk[v]} = \{ 1\}$.  Hence, Bass--Serre theory (Kurosh's theorem, see \cite[Chapter I, Theorem 14]{Se03}) implies that 
\begin{align*}
    \llangle A_T \rrangle = \left(\ast_{g\in X} \llangle A_T \rrangle \cap A_{\st[v]}^g \right) \ast F ,
\end{align*}
where $F$ is free and $X$ is a system of coset representatives of $\llangle A_T \rrangle \backslash A_\Gamma / A_{\st[v]}$.

Since  $A_{\st[v]}= A_T \times A_{\lk[v]}$, we immediately see that 
 $\llangle A_T \rrangle \cap A_{\st[v]} = A_T$. Now using that $\llangle A_T \rrangle$ is normal, it follows that for all $g\in A_\Gamma$, we have $\llangle A_T \rrangle \cap A_{\st[v]}^g= \left(\llangle A_T \rrangle \cap A_{\st[v]}\right)^g= A_T^g$ which finishes our proof.
\qed

\begin{lem}\label{lem-comm-free} 
Let $G=\ast_{i\in I} G_i$ with $G_i$ either free or free abelian. Then its commutator subgroup $[G,G]$ is a free group.
\end{lem}
\Proof Let $q\colon G\to G/[G,G]$ be the abelianization map; we have $\ker(q)=[G,G]$. By Bass--Serre theory, $G$ acts on a tree $T$ so that each vertex stabilizer is  conjugate to $G_i$ for some $i$  and all edge stabilizers are trivial.  This means that for every vertex $V\in T$,  its stabilizer  $G_V$ coincides with $gG_ig^{-1}$ for some $g \in G$ and $i\in I$. Now $\ker(q)$ as a  subgroup of $G$ also acts on $T$. The vertex stablizer of $V$ for the action of $\ker(q)$  is $\ker(q) \cap G_V = (\ker(q) \cap G_i)^g$ since $\ker(q)$ is normal.  As $G_i$ is either free or free abelian, the kernel of $q|_{G_i}$ is free or trivial. Hence  $\ker(q) \cap G_V $ is free (bearing in mind that the trivial group is also free). This implies that $\ker(q)$ acts on the tree $T$ with vertex stabilizers being free groups and edge stabilizers trivial. By Bass--Serre theory, $\ker(q)$ is a free product of free groups and hence is itself free.
\qed

\begin{cor}\label{FJC-normalclosure}
For every standard equivalence class $T=[v]$ of $V(\Gamma)$, the subgroup  $\llangle A_T \rrangle$ and every semidirect product $\llangle A_T \rrangle \rtimes\BZ$ satisfy FJCw.
\end{cor}
\Proof
We begin with $\llangle A_T \rrangle$. Abelianisation gives us a map $p\colon \llangle A_T \rrangle \to \llangle A_T \rrangle^{ab}$. Its image is abelian and hence satisfies FJCw. Its kernel is given by $\left[\llangle A_T \rrangle, \llangle A_T \rrangle\right]$. By Lemma \ref{lem-comm-free}, the commutator subgroup of a free product of free and free abelian groups is a free group and hence satisfies FJCw. Now using Theorem \ref{qut} (4), we need to verify that for every infinite cyclic subgroup $C$ of  $ \llangle A_T \rrangle^{ab}$, the preimage $p^{-1}(C)$ also satisfies FJCw. However, such a preimage is a free-by-cyclic group and thus satisfies FJCw by Theorem A.

The argument for $\llangle A_T \rrangle \rtimes\BZ$ is almost identical: Abelianisation induces a surjective map $q\colon \llangle A_T \rrangle \rtimes\BZ \to \llangle A_T \rrangle^{ab}  \rtimes\BZ$. Its image is solvable, its kernel is again the free group $\left[\llangle A_T \rrangle, \llangle A_T \rrangle \right]$ and for every infinite cyclic subgroup $C$ of $\llangle A_T \rrangle^{ab}  \rtimes\BZ$, the preimage $q^{-1}(C)$ is free-by-cyclic.
\qed

By the work of Laurence \cite{La95} and Servatius \cite{Se89}, the group $\Aut(A_\Gamma)$ is generated by automorphisms of the underlying graph $\Gamma$ together with three kinds of automorphisms which are commonly called transvections, inversions and partial conjugations. Let $\Aut^0(A_\Gamma)$ be the subgroup generated by all the transvections, inversions and partial conjugations of $\Aut(A_\Gamma)$ (but not the graph automorphisms). This is a finite-index subgroup of $\Aut(A_\Gamma)$ and it has the property that whenever $T=[v]$ is a standard equivalence class of $V(\Gamma)$ that is maximal with respect to the standard ordering, every element $f\in \Aut^0(A_\Gamma)$ sends $A_T$ to a conjugate of itself (see \cite[Proposition 3.2]{CV09}).

\textbf{Proof of Theorem C.}
We do an induction on the number of vertices of $\Gamma$.
If $A_\Gamma$ is a free or a free abelian group, $A_\Gamma \rtimes \BZ$ is either free-by-cyclic or solvable and hence satisfies FJCw by Theorem \ref{qut} (6) and (8). If this is not the case, every standard equivalence class forms a proper subset of $V(\Gamma)$. Let $T=[v]$ be a maximal such equivalence class. 

The group $\Aut^0(A_\Gamma)$ has finite index in $\Aut(A_\Gamma)$, so there is a natural number $k$ such that $f^k\in \Aut^0(A_\Gamma)$. The group $A_\Gamma \rtimes_{f^k}\BZ=A_\Gamma \rtimes_f k\BZ$ has finite index in $A_\Gamma \rtimes_f\BZ$, so by Theorem \ref{qut} (6), it satisfies FJCw if and only if $A_\Gamma \rtimes_f\BZ$ does. Hence, we may assume that $f\in \Aut^0(A_\Gamma)$.

As every element of $\Aut^0(A_\Gamma)$ sends $A_T$ to a conjugate of itself, the automorphism $f$ stabilizes $\llangle A_T \rrangle$. Consequently, we get a quotient map $q\colon A_\Gamma \rtimes_f\BZ \to A_\Gamma/\llangle A_T \rrangle \rtimes_{\bar{f}}\BZ$, where $\bar{f}$ is the induced map on $A_\Gamma/\llangle A_T \rrangle$. The group $A_\Gamma/\llangle A_T \rrangle$ is isomorphic to $A_{V(\Gamma) \backslash T}$, so the image of $q$ satisfies FJCw by induction. The kernel of $q$ is $\llangle A_T \rrangle$ and thus satisfies FJCw by Corollary \ref{FJC-normalclosure}. Given an infinite cyclic subgroup $C$ of $A_\Gamma/\llangle A_T \rrangle \rtimes_{\bar{f}}\BZ$, the preimage $q^{-1}(C)$ is either contained in $A_\Gamma$ and hence satisfies FJCw by Theorem \ref{qut} (1) or it is of the form $\llangle A_T \rrangle \rtimes\BZ$ and satisfies FJCw by Corollary \ref{FJC-normalclosure}. The claim now follows again from Theorem \ref{qut} (4).
\qed
\bibliographystyle{amsplain}

\end{document}